\documentclass[11pt]{article}
\usepackage{graphicx,subfigure}
\usepackage{epsfig}
\usepackage{amssymb}
\usepackage[dvipdfm]{hyperref}
\newcommand{\E}{{\cal E}}

\newtheorem{theorem}{Theorem}[section]

\newtheorem{lemma}[theorem]{Lemma}

\def\whitebox{{\hbox{\hskip 1pt
 \vrule height 6pt depth 1.5pt
 \lower 1.5pt\vbox to 7.5pt{\hrule width
    3.2pt\vfill\hrule width 3.2pt}%
 \vrule height 6pt depth 1.5pt
 \hskip 1pt } }}
\def\qed{\ifhmode\allowbreak\else\nobreak\fi\hfill\quad\nobreak
     \whitebox\medbreak}

\newcommand{\ignore}[1]{}

\begin {document}

\baselineskip 16pt
\title{Asymptotic Laplacian-Energy-Like Invariant of Lattices}

\author{\small $  Jia\textrm{-}Bao \ Liu^{1,2},\ \  Xiang\textrm{-}Feng \ Pan^1
 \thanks{Corresonding author. Tel:+0551-6386-1835. \E-mail:~xfpan@ahu.edu.cn(X.-F. Pan),liujiabaoad@163.com(J.-B. Liu),
  hufu@mail.ustc.edu.cn(F.-T. Hu),
 huayanyanhu@163.com(F.-F.Hu).}\ \ ,Fu\textrm{-}Tao \ Hu^{1},\ \  Feng\textrm{-}Feng \ Hu^1
 $\\ \\
\small  $^{1}$ School of Mathematical Sciences, Anhui  University, Hefei 230601, China\\
\small  $^{2}$ Department of Public Courses, Anhui Xinhua
University, Hefei 230088, China\\}

\date{}
\maketitle
\begin{abstract}
Let $\mu_1\ge \mu_2\ge\cdots\ge\mu_n$ denote the Laplacian
eigenvalues of $G$ with $n$ vertices. The Laplacian-energy-like
invariant, denoted by $LEL(G)= \sum_{i=1}^{n-1}\sqrt{\mu_i}$, is a
novel topological index. In this paper, we show that the
Laplacian-energy-like per vertex of various lattices is
independent of the toroidal, cylindrical, and free boundary
conditions. Simultaneously, the explicit asymptotic values of the
Laplacian-energy-like in these lattices are obtained. Moreover,
our approach implies that in general the Laplacian-energy-like per
vertex of other lattices is independent of the boundary
conditions.

\end{abstract}

\section{ Introduction}
 Throughout this paper we concerned with finite undirected
connected simple graphs. Let $G$ be a graph with vertices labelled
$1,2,\dots,n$. The adjacency matrix $A(G)$ of $G$ is an $n\times
n$ matrix with the $(i,j)$-entry equal to 1 if vertices $i$ and
$j$ are adjacent and 0 otherwise. Supposed
$D(G)=diag(d_1(G),d_2(G),\ldots ,d_n(G))$ be the degree diagonal
matrix of $G$, where $d_i(G)$ is the degree of the vertex $i$,
$i=1,2,\ldots,n$. Let $L(G)=D(G)-A(G)$ is called the Laplacian
matrix of $G$. Then, the eigenvalues of $A(G)$ and $L(G)$ are
called eigenvalues and Laplacian eigenvalues of $G$, denoted by
$\lambda_1, \lambda_2,\dots,\lambda_n$  and $\mu_1,
\mu_2,\dots,\mu_n$, respectively. The underlying graph theoretical
definitions and notations  follow {\rm \cite{CZ}.

Gutman has defined the energy of a graph $G$ with $n$ vertices {\rm
\cite{IG}, denoted by $E(G)$, as $$E(G)=
\sum_{i=1}^{n}|\lambda_i|.$$

Let $G$ be a graph of order $n$ with Laplacian spectrum $\mu_1
\geq \mu_2 \geq \cdots \geq \mu_n.$ The Laplacian-energy-like
invariant of G, LEL for short, is defined as
$$  LEL(G)=\sum_{i=1}^{n-1}\sqrt{\mu_i}.$$
The concept of $LEL(G)$ was introduced by J. Liu and B. Liu ({\rm
\cite{LL}, 2008), where it was shown that it has similar features
as graph energy.

Gutman et al. pointed out in {\rm \cite{GZF} that $LEL(G)$ is more
similar to $E(G)$ than to $LE(G)$.
 The Laplacian-energy-like invariant describes well the properties which are accounted by the
majority of molecular descriptors: motor octane number, entropy,
molar volume, molar refraction, particularly the acentric factor
AF parameter, but also more difficult properties like boiling
point, melting point and partition coefficient Log P. {\rm
\cite{LHY,DA}. D. Stevanovi$\acute{c}$ et al. {\rm \cite{DA} have
proved that $LEL(G)$ has the properties as good as the
Randi$\acute{c}$ $\chi$ index and better than the Wiener index
which is a distance based index. Besides, it is well defined
mathematically and shows interesting relationships in particular
classes of graphs, these recommending $LEL(G)$ as a noval and
powerful topological index.

 The index has
attracted extensive attention due to its wide applications in
physics, chemistry, graph theory, etc. {\rm \cite{GM,LY,MA}.
Details on its theory can be found in the survey on the
Laplacian-energy-like invariant {\rm \cite{LHY}, the recent papers
{\rm \cite{DGCZ,DXG}, and the references cited therein. Due to its
structure, lattices are of special interest, especially the
chemical and physical indices of some lattices were studied
extensively, see for instance {\rm \cite{YLZ,DGMS,ZZ}.

One of reasons we investigated $LEL(G)$ of lattices is that some
physical and chemical indices of lattices have been presented in
 {\rm \cite{YLZ,ZZ,LPC}, however, there is seldom results on $LEL(G)$ of
 lattices.
 Another reason is, there are a great deal of analogies between
the properties of $E(G)$ and $LEL(G)$, asymptotic energy of
lattices has investigated in chemical physics{\rm \cite{YZ,YYZ},
it is natural for us to consider asymptotic $LEL(G)$ of various
lattices.

The energy $E(G)$  of toroidal lattices have been studied in {\rm
\cite{YZ}, the Kirchhoff index of some toroidal lattices has also
investigated in {\rm \cite{YLZ}, a very elementary and natural
question is that if $G$ are lattices, how about $LEL(G)$? It is an
interesting problem to study the Laplacian-energy-like of some
lattices with various boundary condition.  Motivated by results
above, we consider the problem of computation of the $LEL(G)$ of
some lattices in this article.

The rest of the paper is organized as follows. In Section 2, we
propose the asymptotic Laplacian-energy-like of square lattices
and give the related explanations. We provide a detailed
derivation of the asymptotic Laplacian-energy-like change due to
edge deletion in Section 3. The asymptotic Laplacian-energy-like
of hexagonal, triangular, and $3^34^2$ lattices with three
boundary conditions are investigated in Section 4. We present
concluding remarks and conclude the paper in Section 5.

\section{Asymptotic Laplacian-energy-like of square lattices}

 Given graphs $G$ and $H$ with vertex sets $U$ and $V$, the cartesian product
 $G\Box H$ of graphs $G$ and $H$ is a graph such that
the vertex set of $G\Box H$ is the cartesian product $U\Box V$;
and any two vertices $(u,u')$ and $(v,v')$ are adjacent in $G\Box
H$ if and only if either $u = v$ and $u'$ is adjacent with $v'$ in
$H$, or $u' = v'$ and $u$ is adjacent with $v$ in $G$ {\rm
\cite{CZ}.

Let $P_m\Box P_n$, $P_m\Box C_n$, and $C_m\Box C_n$, denote the
square lattices with free, cylindrical and toroidal boundary
conditions, respectively, where $P_n$ and $C_n$ denote the path and
the cycle with $n$ vertices. Obviously, $P_m\Box P_n$ is a sequence
of spanning subgraphs of the sequence $P_m\Box C_n$ of finite
graphs, and $P_m\Box C_n$ is a sequence of spanning subgraphs of the
sequence $C_m\Box C_n$ of finite graphs. Particularly,
\begin{eqnarray} \nonumber  &&\lim_{m,n\to \infty}\frac{\left|\{v\in
V(P_m\Box P_n):~d_{P_m\Box P_n}(v) =
d_{C_m\Box C_n}(v)\}\right|}{\mid V(C_m\Box C_n)\mid} \\
&&=\lim_{m,n\to \infty}\frac{\left|\{v\in V(P_m\Box C_n):~d_{P_m\Box
C_n}(v)=d_{C_m\Box C_n}(v)\}\right|}{\mid V(C_m\Box C_n)\mid}
 =1,\nonumber
 \end{eqnarray}
  that is, almost
all vertices of $C_m\Box C_n$ and $P_m\Box C_n$
 (resp. $C_m\Box C_n$ and $P_m\Box P_n$) have the same degrees.
Let $G_1,G_2$ be graphs with adjacency matrices $A_1,$ $A_2$,
degree matrices $D_1,$ $D_2$ and Laplacan matrices $L_1,$¡¡$L_2,$
respectively. Then $L_1=D_1-A_1,L_2=D_2-A_2.$
 If $\mu_i(G_1),\mu_j(G_2)$ are Laplacian eigenvalues of $G_1,G_2$, then
the Laplacian eigenvalues of $G_1\Box G_2$ are all possible sums
$\mu_i(G_1)+\mu_j(G_2)$, as noted in \cite{DC}. On the other hand,
it is well known that the Laplacian eigenvalues of a path $P_m$
and a cycle $C_m$ are $2-2cos\frac{i\pi}{n} (i=0,1,\dots,n-1)$
 and $2-2cos\frac{2\pi j}{n} (j=0,1,\dots,m-1)$
\cite{DC}, respectively. Consequently, the Laplacian eigenvalues
of $P_m\Box P_n$ (resp. $P_m\Box C_n$ and $C_m\Box C_n$) are
$4-2cos\frac{i\pi}{m}-2cos\frac{j\pi}{n},i=0,1,\dots,n-1;j=0,1,\dots,m-1$
(resp.
$4-2cos\frac{i\pi}{m}-2cos\frac{2j\pi}{n},i=0,1,\dots,n-1;j=0,1,\dots,m-1$
and
$4-2cos\frac{2i\pi}{m}-2cos\frac{2j\pi}{n},i=0,1,\dots,n-1;j=0,1,\dots,m-1$).

Therefore, the Laplacian-energy-like per vertex of $P_m\Box P_n$,
$P_m\Box C_n$, and $C_m\Box C_n$ are defined as
$$1. \ \lim_{m\to\infty} \lim_{n\to \infty}\frac{LEL\Big(P_m\Box P_n\Big)}{\mid V(P_m\Box P_n)\mid}
= \lim_{m\to\infty} \lim_{n\to \infty}\frac{1}{mn}
\sum_{i=0}^{m-1}\sum_{j=0}^{n-1}
\sqrt{4-2cos\frac{i\pi}{m}-2cos\frac{j\pi}{n}} $$
$$=\int_0^{1}\int_0^{1} \sqrt{4-2cos\pi x-2cos\pi y}\ ~dxdy
= \frac{1}{\pi^2}\int_0^{\pi}\int_0^{\pi} \sqrt{4-2cos x-2cos y}\
~dxdy\approx 1.9162.$$

$$2. \ \lim_{m\to\infty} \lim_{n\to \infty}\frac{LEL\Big(P_m\Box C_n\Big)}{\mid V(P_m\Box C_n)\mid}
= \lim_{m\to\infty} \lim_{n\to \infty}\frac{1}{mn}
\sum_{i=0}^{m-1}\sum_{j=0}^{n-1}
\sqrt{4-2cos\frac{i\pi}{m}-2cos\frac{2j\pi}{n}} $$
$$=\int_0^{1}\int_0^{1} \sqrt{4-2cos\pi x-2cos2\pi y}\ ~dxdy
= \frac{1}{2\pi^2}\int_0^{\pi}\int_0^{2\pi} \sqrt{4-2cos x-2cos
y}\ ~dxdy \approx1.9162.$$

$$3. \ \lim_{m\to\infty} \lim_{n\to \infty}\frac{LEL\Big(C_m\Box C_n\Big)}{\mid V(C_m\Box C_n)\mid}
= \lim_{m\to\infty} \lim_{n\to \infty}\frac{1}{mn}
\sum_{i=0}^{m-1}\sum_{j=0}^{n-1} \sqrt{4-2cos\frac{2
i\pi}{m}-2cos\frac{2j\pi}{n}} $$
$$=\int_0^{1}\int_0^{1} \sqrt{4-2cos2\pi x-2cos2\pi y}\ ~dxdy
= \frac{1}{4\pi^2}\int_0^{2\pi}\int_0^{2\pi} \sqrt{4-2cosx-2cos
y}\ ~dxdy \approx 1.9162.$$

The numerical integration value in last line is calculated with
software MATLAB calculation.

\vspace{5pt} \noindent {\bf Remark 2.1}\ By using computer software
MATLAB, we can easily get that the above numerical integration
values, which imply that $P_m\Box P_n$, $P_m\Box C_n$ and $C_m\Box
C_n$ have the same asymptotic Laplacian-energy-like $LEL(P_m\Box
P_n)=LEL(P_m\Box C_n)=LEL(C_m\Box C_n)\approx 1.9162mn$ as $m,n$
approach infinity.

\vspace{5pt} \noindent {\bf Remark 2.2}\ The asymptotic
Laplacian-energy-like $LEL(G)$ of square lattices is independent on
the three boundary conditions, i.e., the free, cylindrical and
toroidal boundary conditions.

The phenomenon above is not accidental, the related explanations
are proposed in next section. Our method implies
 that in general the Laplacian-energy-like per vertex
of lattices is independent of the boundary conditions.

\section{ Graph asymptotic Laplacian-energy-like change due to edge deletion}

We propose a detailed derivation of the asymptotic
Laplacian-energy-like change due to edge deletion, in our proof,
some techniques in \cite{YZ} are referred to. Recall some results
which will be used in our discussion. The authors \rm \cite{LL}
verified the following theorem.
\begin{theorem}\label{3-1}(\cite{LL})
 Let $G$ be a graph on $n$ vertices
with $m$ edges, then
$$\sqrt{2m}\leq LEL(G)\leq \sqrt{2}m.$$
The first equality is attained if and only if $G\cong\overline{K_n}$
or $K_2\cup(n-2)K_1,$ and the second equality is attained if and
only if $G\cong rK_2\cup(n-2r)K_1$, where $0\leq r\leq
\lfloor\frac{n}{2}\rfloor.$
\end{theorem}

 According to Theorem 3.1, one can easily get
 that
 \begin{lemma}\label{3-2}
 Let $G$ be a graph on $n$ vertices
with $m$ edges, then $$LEL(G)\leq \sqrt{2}m < 2m.$$
\end{lemma}

The authors \rm \cite{DS,DWS} investigated how the energy of a
graph changes when edges are removed, and obtained the following
consequence.
\begin{lemma}\label{3-3} (\cite{DS,DWS})
 Let $H$ be an induced
subgraph of a graph $G.$ Then
$$E(G)-E(H)\leq E(G-E(H))  \leq E(G)+E(H).$$
\end{lemma}

In addition, the authors \rm \cite{YZ} obtained that
$$\mid E(G)-E(H)\mid \leq E(G-E(H))  \leq E(G)+E(H).$$

 With a similar method, one
can prove the following result.
\begin{lemma}\label{3-4}
 Let $H$ be a subgraph of a graph $G.$
  Then
  $$\mid LEL(G)-LEL(H)\mid \leq LEL(G-E(H))  \leq LEL(G)+LEL(H).$$
\end{lemma}

Consider two graphs $G$ and $H$ $(V(G)\cap V(H)$ may be disjoint),
denoted by $$\Delta(G,H)=\mid E(G)\mid+ \mid E(H)\mid-2\mid
E(G)\cap E(H)\mid,$$ i.e., $\Delta(G,H)$ equals the number of
edges of symmetric difference of $E(G)$ and $E(H)$.

\begin{theorem}\label{3-5}
 Let $\{G_n\}$ and $\{H_n\}$ be two sequences of graphs such that
$$\lim_{n\to \infty}\frac{\Delta(G_n,H_n)}{LEL(G_n)}=0.$$
Then
$$\lim_{n\to \infty}\frac{LEL(H_n)}{LEL(G_n)}=1.$$
\end{theorem}
{\bf Proof.} Let $F_n$ be the subgraph of $G_n$ or $H_n$ induced by
$E(G)\cap E(H).$ Note that
\begin{eqnarray}
\nonumber \left| \frac{LEL(H_n)}{LEL(G_n)}-1 \right|
&&=\left|\frac{LEL(H_n)-LEL(G_n)}{LEL(G_n)}\right|
=\left|\frac{LEL(H_n)-LEL(F_n)+LEL(F_n)-LEL(G_n)}{LEL(G_n)}  \right|\\
&&\nonumber \leq \left|\frac{LEL(G_n)-LEL(F_n)}{LEL(G_n)}\right|
+\left|\frac{LEL(H_n)-LEL(F_n)}{LEL(G_n)}\right|.\\ \nonumber
\end{eqnarray}
By Lemma 3.2 and Lemma 3.4,
$$\mid LEL(G_n)-LEL(F_n)\mid \leq LEL(G_n-E(F_n))  \leq 2\left|E(G_n)\right|-2\left|E(F_n)\right|,$$
$$\mid LEL(H_n)-LEL(F_n)\mid \leq LEL(H_n-E(F_n))  \leq 2\left|E(H_n)\right|-2\left|E(F_n)\right|.$$
Consequently,
$$\left| \frac{LEL(H_n)}{LEL(G_n)}-1 \right|\leq \frac{2\Big(\mid E(G_n)\mid+ \mid E(H_n)\mid-2\mid
E(F_n)\mid\Big)  }{LEL(G_n)} =\frac{2\Delta(G_n,H_n)
}{LEL(G_n)},$$ implying the theorem holds. \hfill $\blacksquare$

Based on Theorem 3.5, one can straightforwardly arrive to that

\begin{theorem}\label{3-6}
 Let $\{G_n\}$ be a sequence of finite simple graphs with bounded average degree such that
$$\lim_{n\to \infty}\left|V(G_n)\right|=\infty, \lim_{n\to
\infty}\frac{LEL(G_n)}{\left|V(G_n)\right|}=h\neq0.$$ Let $\{H_n\}$
be a sequence of spanning subgraphs of $\{G_n\}$ such that
$$\lim_{n\to\infty}\frac{\left|v \in V(H_n):d_{H_n(v)}=d_{G_n(v)}\right|}{\left|V(G_n)\right|}=1,$$
 then
$$\lim_{n\to \infty}\frac{LEL(H_n)}{\left|V(G_n)\right|}=h.$$
That is, $G_n$ and $H_n$ have the same asymptotic
Laplacian-energy-like.
\end{theorem}

A direct sequence of Theorem 3.6 is that $P_n\square
P_n$,$P_n\square C_n$, and $C_n\square C_n$ have the same asymptotic
Laplacian-energy-like which is shown in the introduction. More
generally, by Theorem 3.6, we have

\vspace{5pt} \noindent {\bf Remark 3.7}\ Let $G_i=P_n$ or
$G_i=C_n$, $i=1,2,\dots,k,$ and $k$ be a constant. If $n$ is
sufficiently large, then the asymptotic Laplacian-energy-like of
the $n$-dimensional lattices
$$LEL(G_1\square G_2 \square  \dots\square G_k)
\approx \frac{n^k}{\pi^k}\int_0^{\pi}\int_0^{\pi}\dots
\int_0^{\pi}\sqrt{\sum_{i=1}^k\Big(2-2cosx_i\Big)} ~dx_1~dx_2\dots
dx_k.$$

\vspace{5pt} \noindent {\bf Remark 3.8}\   Theorem 3.6 provides a
very effective approach to handle the asymptotic
Laplacian-energy-like of a graph with bounded average degree.

 We will use the approach above to
deal with the asymptotic Laplacian-energy-like of some lattices in
the next subsections.

\section{Asymptotic Laplacian-energy-like of some lattices}

\subsection{The hexagonal lattice}

Our notation for the hexagonal lattices follows {\rm
\cite{YZ,SW}}. The hexagonal lattices with toroidal, cylindrical
and free boundary conditions, denoted by $H^t(m,n)$,  $H^c(m,n)$,
and $H^f(m,n)$ are illustrated in Figure 1, respectively, where
$(a_1, b_1), (a_2, b_2), \dots, (a_{m+1}, b_{m+1}), (a_1, c_1^*),
(c_1, c_2^*),(c_2, c_3^*), \dots, (c_{n-1}, c_n^*),\\ (c_n,
b_{m+1})$ are edges in $H^t(m,n)$, and $(a_1, b_1), (a_2, b_2),
\dots, (a_{m+1}, b_{m+1})$ are edges in $H^c(m,n)$(see Figure
1(b)). If we delete edges $(a_1, b_1), (a_2, b_2), \dots,
(a_{m+1}, b_{m+1})$ from $H^c(m,n)$, then the hexagonal lattice,
denoted by $H^f(m,n)$, with free boundary condition is obtained
(see Figure 1(c)).

\begin{figure}[ht]
\centering
  \includegraphics[width=\textwidth]{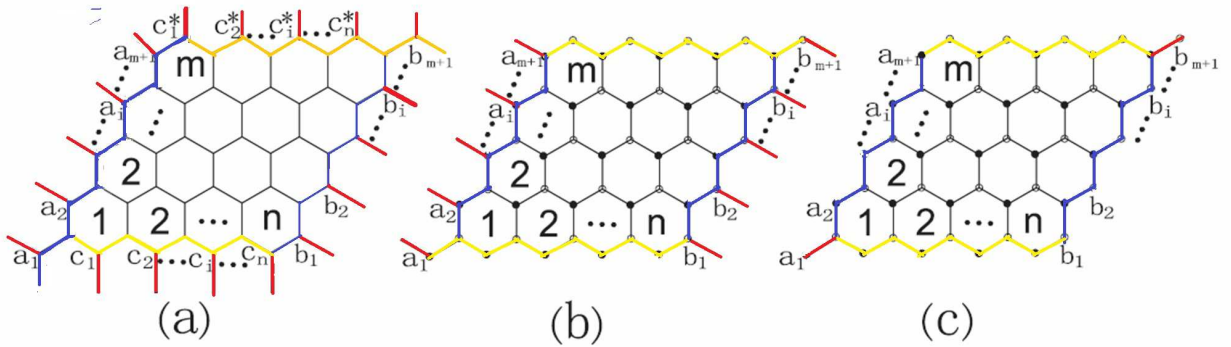}
\caption{(a) The hexagonal lattice $H^t(m,n)$ with toroidal
boundary condition; (b) the hexagonal lattice $H^c(m,n)$ with
cylindrical boundary condition; (c) the hexagonal lattice
$H^f(m,n)$ with free boundary condition.}
\end{figure}

A recent approach to compute the Laplacian eigenvalues of
$H^{t}(m,n)$ can be found in {\rm \cite{GM,DGMS,JS}. Using
Equation (6.2.2) in {\rm \cite{SW}, the Laplacian matrix
$L(H^{t}(m,n))$ of $H^{t}(m,n)$ is similar to the block diagonal
matrix whose diagonal blocks are
$$B_{ij}=\left(%
\begin{array}{cc}
  3 & -1-\omega_{n+1}^{-i}- \omega_{m+1}^{j}\\
 -1-\omega_{n+1}^{i}- \omega_{m+1}^{-j} & 3  \\
 \end{array}%
\right)$$
 where $\omega_s=cos\frac{2\pi}{s}+i
sin\frac{2\pi}{s},i=0,1,\dots,m;j=0,1,\dots,n.$ Hence the
Laplacian eigenvalues of $H^{t}(m,n)$ are
$$3\pm \sqrt{3+2cos\frac{2\pi i}{m+1}
+2cos\frac{2\pi j}{n+1}+2cos\Big(\frac{2\pi i}{m+1}+\frac{2\pi
j}{n+1}\Big)}, \  \ i=0,1,\dots,m;j=0,1,\dots,n.$$

By the definition of the Laplacian-energy-like, it is not difficult
to prove the following result:

\begin{theorem}\label{4-1}
For the hexagonal lattices $H^t(m,n), H^c(m,n)$ and $H^f(m,n)$
with toroidal, cylindrical, and free boundary conditions. Then
 $$1.~~\lim_{m\to\infty}
\lim_{n\to\infty}\frac{
 LEL\Big(H^t(m,n)\Big)}{2(m+1)(n+1)} =\lim_{m\to\infty}
\lim_{n\to\infty}\frac{ LEL\Big(H^c(m,n)\Big)}{2(m+1)(n+1)}
=\lim_{m\to\infty} \lim_{n\to\infty}\frac{
 LEL\Big(H^f(m,n)\Big)}{2(m+1)(n+1)}\approx  1.6437.$$\\
  $2.~LEL\Big(H^{t}(m,n)\Big)=LEL\Big(H^{c}(m,n)\Big)=LEL\Big(H^{f}(m,n)\Big)  \approx
3.2714(m+1)(n+1).$
\end{theorem}

{\bf Proof.} By definitions of $H^t(m,n), H^c(m,n)$ and
$H^f(m,n)$, it is obvious that $H^c(m,n)$ and $H^f(m,n)$ are
spanning subgraphs of $H^t(m,n)$. Furthermore, the degree of
almost all vertices of $H^t(m,n), H^c(m,n)$ and $H^f(m,n)$ are 3.
Hence, by Theorem 3.6, one can obtian that
$$\lim_{m\to \infty}\lim_{n\to \infty}\frac{
 LEL\Big(H^t(m,n)\Big)}{2(m+1)(n+1)} =\lim_{m\to \infty}\lim_{n\to
\infty}\frac{ LEL\Big(H^c(m,n)\Big)}{2(m+1)(n+1)} =\lim_{m\to
\infty}\lim_{n\to \infty}\frac{
 LEL\Big(H^f(m,n)\Big)}{2(m+1)(n+1)}.$$

It suffices to prove that
\begin{eqnarray}
&&\nonumber\lim_{m\to\infty} \lim_{n\to
\infty}\frac{ LEL\Big(H^{t}(m,n)\Big)}{2(m+1)(n+1)}\\
  \nonumber &=& \lim_{m\to\infty} \lim_{n\to \infty}
  \frac{1}{2(m+1)(n+1)}\sum_{i=0}^{m}\sum_{j=0}^{n}
\sqrt{ 3+ \sqrt{3+2cos\alpha_i +2cos\beta_j+2cos(\alpha_i+\beta_j)} }\\
  \nonumber & &+\lim_{m\to\infty} \lim_{n\to \infty}\frac{1}{2(m+1)(n+1)}
  \sum_{i=0}^{m}\sum_{j=0}^{n} \sqrt{ 3- \sqrt{3+2cos\alpha_i +2cos\beta_j+2cos(\alpha_i+\beta_j)} }\\
  \nonumber &=&\frac{1}{2}\int_{0}^{1}\int_{0}^{1}
\sqrt{ 3+ \sqrt{3+2cos 2\pi x+2cos 2\pi y+2cos 2\pi(x+y)} } \ ~dxdy \\
  \nonumber & &+\frac{1}{2}\int_{0}^{1}\int_{0}^{1}
\sqrt{ 3- \sqrt{3+2cos 2\pi x+2cos 2\pi y+2cos 2\pi(x+y)} } \ ~dxdy\\
 \nonumber &=&\frac{1}{8\pi^2}\int_{0}^{2\pi}\int_{0}^{2\pi}
\sqrt{ 3+ \sqrt{3+2cos x+2cos y+2cos(x+y)} } \ ~dxdy \\
  \nonumber & &+\frac{1}{8\pi^2}\int_{0}^{2\pi}\int_{0}^{2\pi}
\sqrt{ 3- \sqrt{3+2cos x+2cos y+2cos(x+y)} } \ ~dxdy\\
\nonumber &\approx & 1.6437 .
\end{eqnarray}

The above numerical integration value implies that the hexagonal
lattices $H^t(m,n), H^c(m,n)$ and $H^f(m,n)$ with toroidal,
cylindrical, and free boundary conditions have the same asymptotic
Laplacian-energy-like, i.e.,
$LEL\Big(H^{t}(m,n)\Big)=LEL\Big(H^{c}(m,n)\Big)=LEL\Big(H^{f}(m,n)\Big)
\approx 3.2714(m+1)(n+1)$ as $m,n$ tends to infinity. \hfill
$\blacksquare$

\subsection{The 3.12.12 lattice}
Our notation for the 3.12.12 lattices follows {\rm \cite{LY,ZZ}}.
The 3.12.12 lattice with toroidal boundary condition, denoted by
$J^{t}(m,n)$, is the graph illustrated in Figure 2(a).

\begin{figure}[h]

\subfigure[]{
\begin{minipage}[t]{0.32\linewidth}
\centering
 \includegraphics[width=\textwidth]{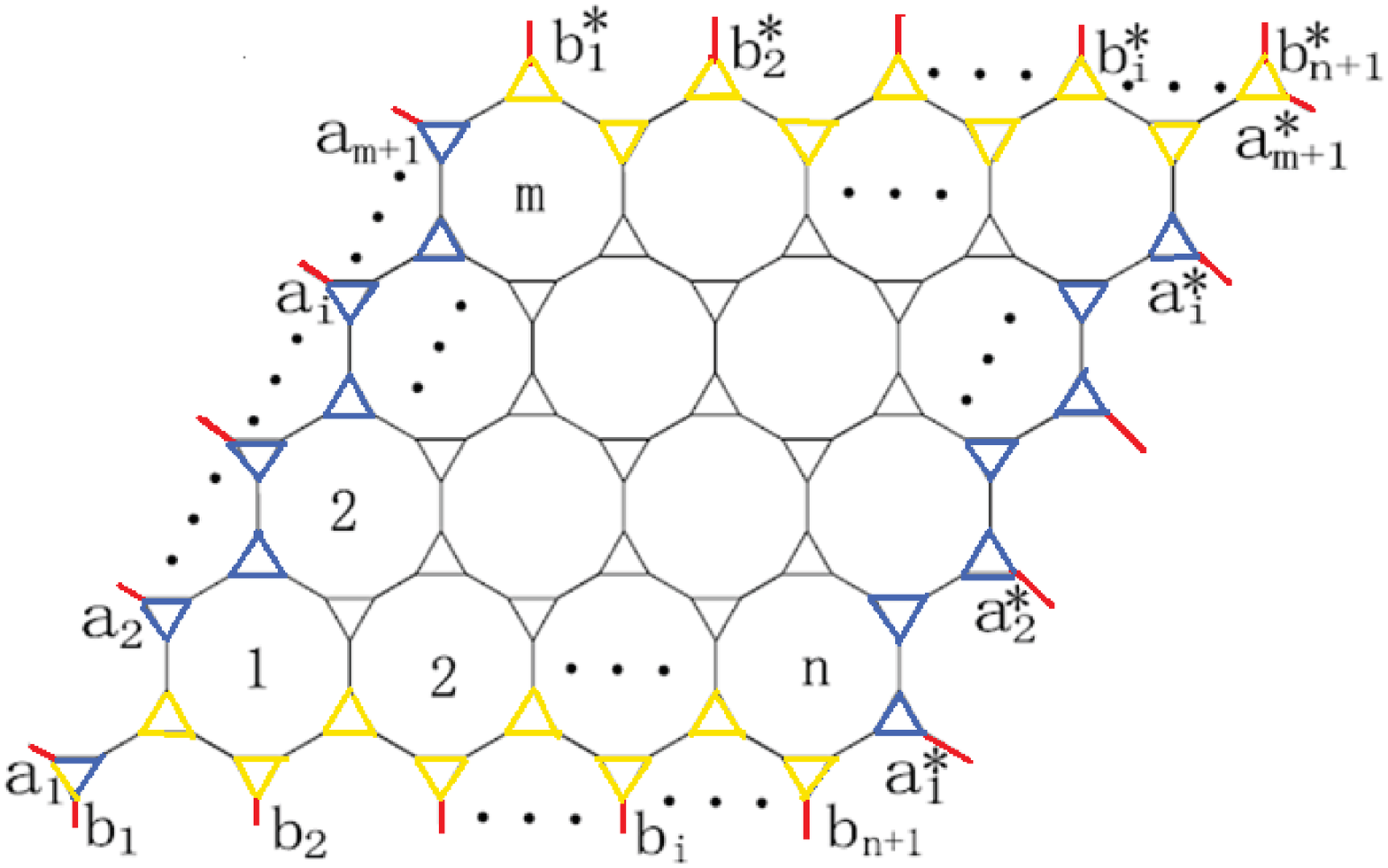}
\end{minipage}}%
\subfigure[]{
\begin{minipage}[t]{0.32\linewidth}
\centering
 \includegraphics[width=\textwidth]{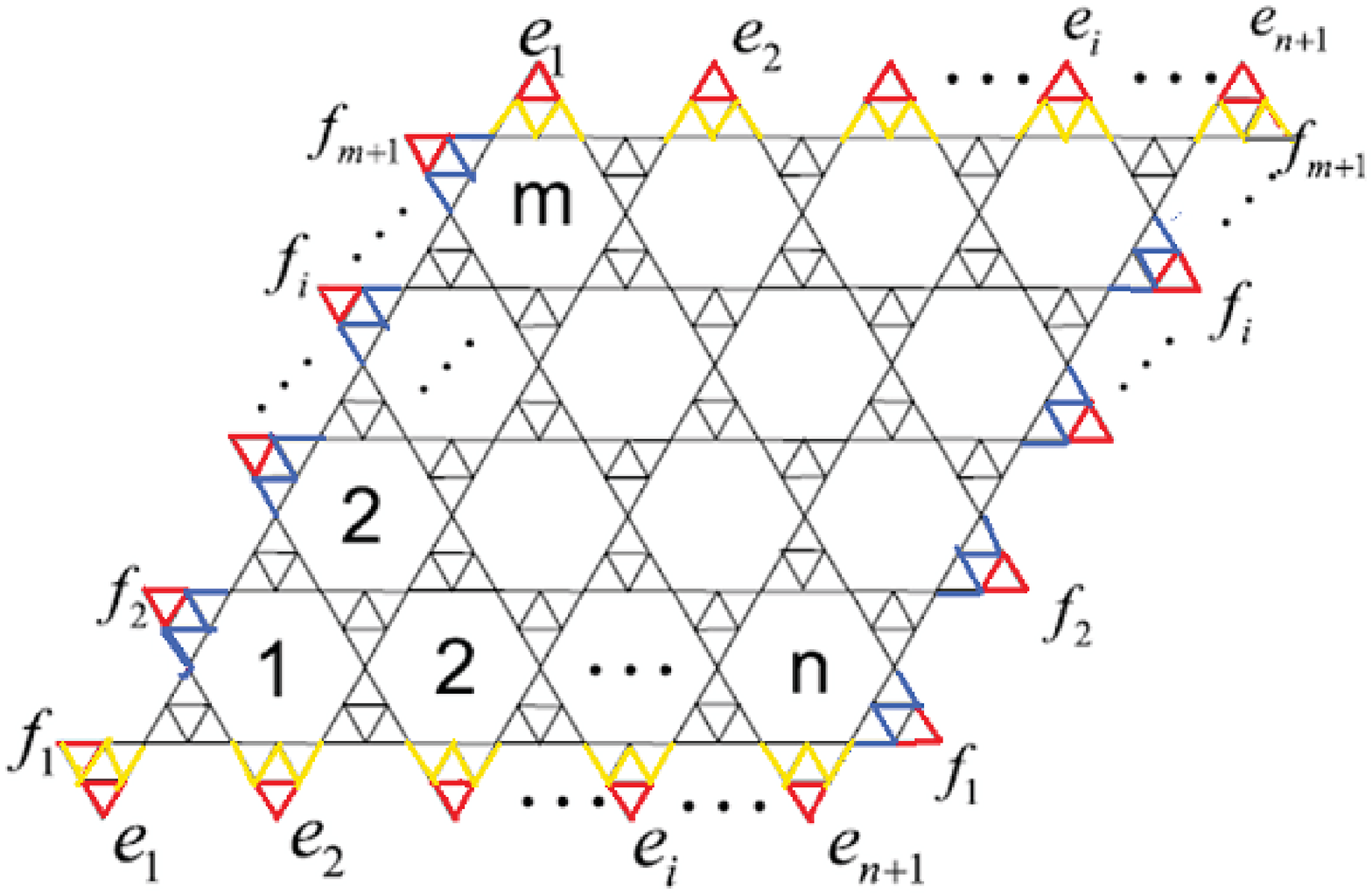}
\end{minipage}}%
\subfigure[]{
\begin{minipage}[t]{0.32\linewidth}
\centering
 \includegraphics[width=0.8\textwidth]{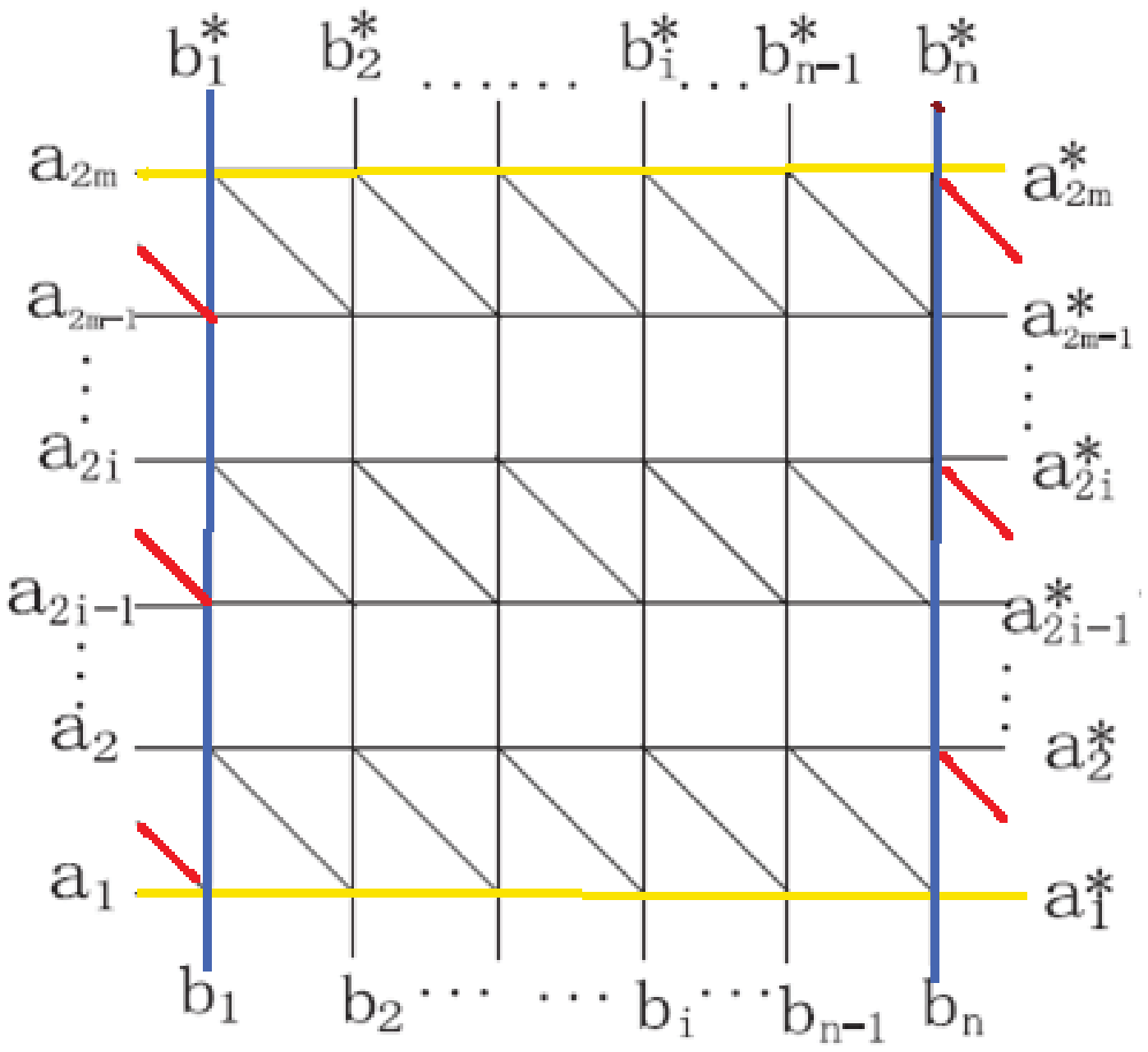}
\end{minipage}}
\caption{(a) The 3-12-12 lattice with toroidal boundary condition;
(b) The triangular kagom$\acute{e}$ lattice with toroidal boundary
condition; (c) the $3^3.4^2$ lattice with toroidal boundary
condition.}
\end{figure}

Recently, the adjacency spectrum of 3.12.12 lattice has been
proposed in {\rm \cite{LY}.

\begin{theorem}\label{2-2} {\rm \cite{LY}}
 Let $J^{t}(m,n)$ be the 3.12.12 lattice with toroidal boundary condition. Then the adjacency spectrum is
\begin{eqnarray}
\nonumber Spec_A(J^{t}(m,n))&=&\textbf{\{}\underbrace{-2,-2,\dots,-2}_{(m+1)(n+1)},\underbrace{0,0,\dots,0}_{(m+1)(n+1)}\textbf{\}}\\
\nonumber &&  \bigcup\left\{ \frac{  1\pm \sqrt
{13\pm4\sqrt{3+2cos\alpha_i+2cos\beta_j+2cos\Big(\alpha_i+\beta_j\Big)}}
    }{2}:0\le i\le m, 0\le j\le n  \right\},\\ \nonumber
    \end{eqnarray}
    where $\alpha_i=\frac{2\pi i}{m+1}, \beta_j=\frac{2\pi
j}{n+1},  i=0,1,\dots,m;j=0,1,\dots,n.$
\end{theorem}

The following result has been obtained in {\rm \cite{NLB}, which
is an important relationship between $Spec_A(G)$ and $Spec_L(G)$.
Suppose that $G$ is an $r$-regular graph with $n$ vertices and
$Spec_A(G)=\{\lambda_1, \lambda_2, \dots, \lambda_n\}.$ Then
$$ Spec_L(G)=\Big\{r-\lambda_1, r-\lambda_2, \dots, r-\lambda_n\Big\}.$$

Note that $J^{t}(m,n)$ is the line graph of the subdivision of
$H^t (m,n)$ which is a $3$-regular graph with $2(m + 1)(n + 1)$
vertices and $J^{t}(m,n)$ has $6(m + 1)(n + 1)$ vertices. Hence,
we arrive to that

\begin{theorem}\label{4-3}
 Let $J^{t}(m,n)$ be the 3.12.12 lattice with toroidal boundary condition. Then the Laplacian spectrum is
\begin{eqnarray}
\nonumber Spec_L(J^{t}(m,n))&=&\textbf{\{}\underbrace{5,5,\dots,5}_{(m+1)(n+1)},\underbrace{3,3,\dots,3}_{(m+1)(n+1)}\textbf{\}}\\
\nonumber &&  \bigcup\left\{\frac{  5\pm \sqrt
{13\pm4\sqrt{3+2cos\alpha_i
+2cos\beta_j+2cos\Big(\alpha_i+\beta_j\Big)}}
    }{2}:0\le i\le m, 0\le j\le n \right\},\\ \nonumber
    \end{eqnarray}
    where $\alpha_i=\frac{2\pi i}{m+1}, \beta_j=\frac{2\pi
j}{n+1},  i=0,1,\dots,m;j=0,1,\dots,n.$
\end{theorem}

By the definition of the Laplacian-energy-like $LEL(G)$, one can
easily arrive to the following theorem.
\begin{theorem}\label{4-4}
For $J^{t}(m,n), J^c(m,n)$ and $J^f(m,n)$ with toroidal,
cylindrical, and free boundary conditions. Then
$$1.~~ \lim_{m\to\infty} \lim_{n\to
\infty}\frac{LEL\Big(J^{t}(m,n)\Big)}{6(m+1)(n+1)}=
 \lim_{m\to\infty}
\lim_{n\to\infty}\frac{LEL\Big(J^{c}(m,n)\Big)}{6(m+1)(n+1)}=
\lim_{m\to \infty}
\lim_{n\to\infty}\frac{LEL\Big(J^{f}(m,n)\Big)}{6(m+1)(n+1)}\approx
 1.3375 .$$\\
 $2.~~LEL\Big(J^t
(m,n)\Big)= LEL\Big(J^c (m,n)\Big)=LEL\Big(J^f (m,n)\Big)\approx
8.0250(m+1)(n+1).$
\end{theorem}
{\bf Proof.} By definitions of $H^t(m,n), H^c(m,n)$ and
$H^f(m,n)$, one can know that $H^c(m,n)$ and $H^f(m,n)$ are
spanning subgraphs of $H^t(m,n)$. Furthermore, the degree of
almost all vertices of $H^t(m,n), H^c(m,n)$ and $H^f(m,n)$ are 3.
Therefore, by Theorem 3.6 it is not difficult to arrive to that
$$ \lim_{m\to\infty} \lim_{n\to
\infty}\frac{LEL\Big(J^{t}(m,n)\Big)}{6(m+1)(n+1)}=
 \lim_{m\to \infty}
\lim_{n\to\infty}\frac{LEL\Big(J^{c}(m,n)\Big)}{6(m+1)(n+1)}=
\lim_{m\to \infty}
\lim_{n\to\infty}\frac{LEL\Big(J^{f}(m,n)\Big)}{6(m+1)(n+1)}.$$

It suffices to prove that
\begin{eqnarray}
\nonumber && \lim_{m\to\infty} \lim_{n\to
\infty}\frac{LEL\Big(J^{t}(m,n)\Big)}{6(m+1)(n+1)}\\
\nonumber &=&\lim_{m\to\infty} \lim_{n\to
\infty}\frac{1}{12(m+1)(n+1)}\sum_{i=0}^{m}\sum_{j=0}^{n}
\sqrt{5- \sqrt {13-4\sqrt{3+2cos\alpha_i
+2cos\beta_j+2cos\Big(\alpha_i+\beta_j\Big)}
    }}\\  \nonumber\\
     \nonumber  &&+\lim_{m\to\infty} \lim_{n\to
\infty}\frac{1}{12(m+1)(n+1)}\sum_{i=0}^{m}\sum_{j=0}^{n}
\sqrt{5- \sqrt {13+4\sqrt{3+2cos\alpha_i
+2cos\beta_j+2cos\Big(\alpha_i+\beta_j\Big)}
    }}\\  \nonumber
      \nonumber  &&+\lim_{m\to\infty} \lim_{n\to
\infty}\frac{1}{12(m+1)(n+1)}\sum_{i=0}^{m}\sum_{j=0}^{n}
\sqrt{5+ \sqrt {13-4\sqrt{3+2cos\alpha_i
+2cos\beta_j+2cos\Big(\alpha_i+\beta_j\Big)}
    }}\\  \nonumber
      \nonumber  &&+\lim_{m\to\infty} \lim_{n\to
\infty}\frac{1}{12(m+1)(n+1)}\sum_{i=0}^{m}\sum_{j=0}^{n}
\sqrt{5+ \sqrt {13+4\sqrt{3+2cos\alpha_i
+2cos\beta_j+2cos\Big(\alpha_i+\beta_j\Big)}
    }}\\ \nonumber
    &&+\frac{\sqrt{3}+\sqrt{5}}{6}\\  \nonumber
    \nonumber &=&\frac{1}{12}\int_{0}^{1}\int_{0}^{1}\sqrt{5-
\sqrt {13-4\sqrt{3+2cos 2\pi x +2cos 2\pi y+2cos 2\pi \Big(x+y\Big)}
    }}~dxdy\\  \nonumber\\
     \nonumber  &&+\frac{1}{12}\int_{0}^{1}\int_{0}^{1}\sqrt{5-
\sqrt {13+4\sqrt{3+2cos 2\pi x +2cos 2\pi y+2cos 2\pi \Big(x+y\Big)}
    }}~dxdy\\  \nonumber
      \nonumber  &&+\frac{1}{12}\int_{0}^{1}\int_{0}^{1}\sqrt{5+
\sqrt {13-4\sqrt{3+2cos 2\pi x +2cos 2\pi y+2cos 2\pi \Big(x+y\Big)}
    }}~dxdy\\  \nonumber
      \nonumber  &&+\frac{1}{12}\int_{0}^{1}\int_{0}^{1}\sqrt{5+
\sqrt {13+4\sqrt{3+2cos 2\pi x +2cos 2\pi y+2cos 2\pi \Big(x+y\Big)}
    }}~dxdy+\frac{\sqrt{3}+\sqrt{5}}{6} \\  \nonumber
   \nonumber &\approx & 1.3375 .
\end{eqnarray}

The above numerical integration value implies that
$J^{t}(m,n),J^{c}(m,n)$ and $J^{f}(m,n)$ have the same asymptotic
Laplacian-energy-like, i.e., $LEL\Big(J^t (m,n)\Big)= LEL\Big(J^c
(m,n)\Big)=LEL\Big(J^f (m,n)\Big)\\\approx 8.0250(m+1)(n+1)$ as
$m,n$ tends to infinity. \hfill $\blacksquare$

\subsection{The triangular kagom$\acute{e}$ lattice}

The triangular kagom$\acute{e}$ lattice with toroidal boundary
condition, denoted by $TKL^t(m,n)$, is depicted in Figure 2(b).
Ising spins and XXZ/Ising spins on the $TKL^t(m,n)$ have been
studied in {\rm \cite{YLC,SLM}}. In order to obtain the
Laplacian-energy-like of the The triangular kagom$\acute{e}$
lattice, we recall the spectrum and the Laplacian spectrum of
$TKL^{t}(m,n)$.

\begin{theorem}\label{2-3}{\rm
\cite{LY}}
 Let $\alpha_i=\frac{2\pi i}{m+1}, \beta_j=\frac{2\pi
j}{n+1},  i=0,1,\dots,m;j=0,1,\dots,n.$ Then the spectrum and the
Laplacian spectrum of $TKL^{t}(m,n)$ are
\begin{eqnarray}
\nonumber
Spec_A(TKL^{t}(m,n))&=&\textbf{\{}\underbrace{-2,-2,\dots,-2}_{3(m+1)(n+1)},
\underbrace{-1,-1,\dots,-1}_{(m+1)(n+1)},\underbrace{1,1,\dots,1}_{(m+1)(n+1)}\textbf{\}}\\
\nonumber &&  \bigcup \left\{ \frac{ 3\pm \sqrt
{13\pm4\sqrt{3+2cos\alpha_i
+2cos\beta_j+2cos\Big(\alpha_i+\beta_j\Big)}}
    }{2}:0\le i\le m, 0\le j\le n  \right\},\\ \nonumber
    \end{eqnarray}
    and
\begin{eqnarray}
\nonumber
Spec_L(TKL^{t}(m,n))&=&\textbf{\{}\underbrace{6,6,\dots,6}_{3(m+1)(n+1)},
\underbrace{3,3,\dots,3}_{(m+1)(n+1)},\underbrace{5,5,\dots,5}_{(m+1)(n+1)}\textbf{\}}\\
\nonumber &&  \bigcup\left\{ \frac{ 5\pm \sqrt
{13\pm4\sqrt{3+2cos\alpha_i
+2cos\beta_j+2cos\Big(\alpha_i+\beta_j\Big)}}
    }{2}:0\le i\le m, 0\le j\le n  \right\}.\\ \nonumber
\end{eqnarray}
\end{theorem}

 Note that the triangular $kagom\acute{e}$ lattice is the line graph
of the $3.12.12$ lattice and $TKL^{t}(m,n)$ is a $4$-regular graph
with $9(m + 1)(n + 1)$ vertices.

Consequently,  we can easily get the following Theorem.

\begin{theorem}\label{3-5}
 For $TKL^{t}(m,n), TKL^{c}(m,n)$ and $TKL^f(m,n)$ with toroidal,
cylindrical, and free boundary conditions. Then
$$1.\lim_{m\to\infty} \lim_{n\to
\infty}\frac{LEL(TKL^{t}(m,n))}{9(m+1)(n+1)}= \lim_{m\to\infty}
\lim_{n\to \infty}\frac{LEL(TKL^{c}(m,n))}{9(m+1)(n+1)}=
\lim_{m\to\infty} \lim_{n\to
\infty}\frac{LEL(TKL^{f}(m,n))}{9(m+1)(n+1)}\approx 1.7082.$$\\
$2.~LEL\Big(TKL^{t}(m,n)\Big)=LEL\Big(TKL^{c}(m,n)\Big)=LEL\Big(TKL^{f}(m,n)\Big)\approx
15.3738(m+1)(n+1).$
\end{theorem}

{\bf Proof.} By definitions of $TKL^{t}(m,n), TKL^{c}(m,n)$ and
$TKL^f(m,n)$, one can know that $TKL^{c}(m,n)$ and $TKL^f(m,n)$
are spanning subgraphs of $TKL^{t}(m,n)$. Furthermore, the degree
of almost all vertices of $TKL^{t}(m,n), TKL^{c}(m,n)$ and
$TKL^f(m,n)$ are 4. Therefore, by Theorem 3.6 it is not difficult
to arrive to that
$$\lim_{m\to\infty} \lim_{n\to
\infty}\frac{LEL\Big(TKL^{t}(m,n)\Big)}{9(m+1)(n+1)}=
\lim_{m\to\infty} \lim_{n\to
\infty}\frac{LEL\Big(TKL^{c}(m,n)\Big)}{9(m+1)(n+1)}=
\lim_{m\to\infty} \lim_{n\to
\infty}\frac{LEL\Big(TKL^{f}(m,n)\Big)}{9(m+1)(n+1)}.$$
 It suffices to prove that
\begin{eqnarray}
\nonumber && \lim_{m\to\infty} \lim_{n\to
\infty}\frac{LEL\Big(TKL^{t}(m,n)\Big)}{9(m+1)(n+1)}\\
\nonumber &=&\lim_{m\to\infty} \lim_{n\to
\infty}\frac{1}{18(m+1)(n+1)}\sum_{i=0}^{m}\sum_{j=0}^{n} \sqrt{5-
\sqrt {13-4\sqrt{3+2cos\alpha_i
+2cos\beta_j+2cos\Big(\alpha_i+\beta_j\Big)}
    }}\\  \nonumber\\
     \nonumber  &&+\lim_{m\to\infty} \lim_{n\to
\infty}\frac{1}{18(m+1)(n+1)}\sum_{i=0}^{m}\sum_{j=0}^{n}
\sqrt{5- \sqrt {13+4\sqrt{3+2cos\alpha_i
+2cos\beta_j+2cos\Big(\alpha_i+\beta_j\Big)}
    }}\\  \nonumber
      \nonumber  &&+\lim_{m\to\infty} \lim_{n\to
\infty}\frac{1}{18(m+1)(n+1)}\sum_{i=0}^{m}\sum_{j=0}^{n}
\sqrt{5+ \sqrt {13-4\sqrt{3+2cos\alpha_i
+2cos\beta_j+2cos\Big(\alpha_i+\beta_j\Big)}
    }}\\  \nonumber
      \nonumber  &&+\lim_{m\to\infty} \lim_{n\to
\infty}\frac{1}{18(m+1)(n+1)}\sum_{i=0}^{m}\sum_{j=0}^{n}
\sqrt{5+ \sqrt {13+4\sqrt{3+2cos\alpha_i
+2cos\beta_j+2cos\Big(\alpha_i+\beta_j\Big)}
    }}\\  \nonumber
    &&+\frac{\sqrt{3}+\sqrt{5}}{9}+\frac{\sqrt{6}}{3}\\  \nonumber
    \nonumber &=&\frac{1}{18}\int_{0}^{1}\int_{0}^{1}\sqrt{5-
\sqrt {13-4\sqrt{3+2cos 2\pi x +2cos 2\pi y+2cos 2\pi \Big(x+y\Big)}
    }}~dxdy\\  \nonumber\\
     \nonumber  &&+\frac{1}{18}\int_{0}^{1}\int_{0}^{1}\sqrt{5-
\sqrt {13+4\sqrt{3+2cos 2\pi x +2cos 2\pi y+2cos 2\pi \Big(x+y\Big)}
    }}~dxdy\\  \nonumber
      \nonumber  &&+\frac{1}{18}\int_{0}^{1}\int_{0}^{1}\sqrt{5+
\sqrt {13-4\sqrt{3+2cos 2\pi x +2cos 2\pi y+2cos 2\pi \Big(x+y\Big)}
    }}~dxdy\\  \nonumber
      \nonumber  &&+\frac{1}{18}\int_{0}^{1}\int_{0}^{1}\sqrt{5+
\sqrt {13+4\sqrt{3+2cos 2\pi x +2cos 2\pi y+2cos 2\pi \Big(x+y\Big)}
    }}~dxdy+\frac{\sqrt{3}+\sqrt{5}}{9}+\frac{\sqrt{6}}{3} \\  \nonumber
   \nonumber &\approx &  1.7082.
\end{eqnarray}

The above numerical integration value implies that $TKL^{t}(m,n)$,
$TKL^{c}(m,n)$ and $TKL^{f}(m,n)$ have the same asymptotic
Laplacian-energy-like, i.e.,
$LEL\Big(TKL^{t}(m,n)\Big)=LEL\Big(TKL^{c}(m,n)\Big)\\=LEL\Big(TKL^{f}(m,n)\Big)\approx
15.3738(m+1)(n+1)$  as $m,n$ tends to infinity. \hfill
$\blacksquare$

\subsection{The $3^3.4^2$ lattice}

The $3^3.4^2$ lattice with toroidal boundary condition, denoted by
$M^{t}(n,2m)$, can be constructed by starting with a $2m\times n$
square lattice and adding a diagonal edge connecting the vertices,
i.e., the upper left to the lower right corners of each square in
every other row as shown in Figure 2(c), where $a_1 = b_1, a_{2m}
= b_1^*, a_1^* = b_n, a_{2m}^* = b_n^*,$ and $(a_1, a_1^* ), (a_2,
a_2^*), \dots, (a_{2m}, a_{2m}^*), (b_1, b_1^*), (b_2, b_2^*),
\dots, (b_n, b_n^*),\\ (a_1, a_2^*), (a_3, a_4^* ),$ $ \dots,
(a_{2m-2}, a_{2m}^*)$ are edges in $M^t(n,2m)$.

Let $A(C_{2m})$ be the adjacency matrix of cycle $C_{2m}$, using
the result in {\rm \cite{YZ},the adjacency matrix $A(M^{t}(n,2m))$
of $M^{t}(n,2m)$ has the following form by a suitable labelling of
vertices of $M^{t}(n,2m)$:$$
A\Big(M^{t}(n,2m)\Big)=\left(%
\begin{array}{cccccc}
  A(C_{2m}) & I_{2m}+F_{2m} & 0 &               \dots & 0 & I_{2m}+F_{2m}^T \\
  I_{2m}+F_{2m}^T & A(C_{2m}) & I_{2m}+F_{2m} & \dots & 0 & 0 \\
  0 & I_{2m}+F_{2m}^T & A(C_{2m}) &             \dots & 0 & 0 \\
  \vdots & \vdots & \ddots & \ddots & \ddots & \vdots \\
  0 & 0 & 0 & \dots & A(C_{2m}) & I_{2m}+F_{2m} \\
  I_{2m}+F_{2m} & 0 & 0 & \dots & I_{2m}+F_{2m}^T & A(C_{2m}) \\
\end{array}%
\right)_{n\times n}$$

Notice that $M^{t}(n,2m)$ is an $r$-regular graph. Let
$L(M^{t}(n,2m))$ be the Laplacian matrix of $M^{t}(n,2m)$, it is
not difficult to obtain that $L(M^{t}(n,2m))$ is similar to the
block diagonal matrix whose diagonal blocks are
$$L_{ij}=\left(%
\begin{array}{cc}
  5-\omega_n^i -\omega_n^{-i}  &  -1-\omega_n^{-i} -\omega_m^{j} \\
   -1-\omega_n^i -\omega_m^{-j} &  5-\omega_n^i -\omega_n^{-i} \\
\end{array}%
\right)$$
 where $\omega_s=cos\frac{2\pi}{s}+i
sin\frac{2\pi}{s},i=0,1,\dots,m-1;j=0,1,\dots,n-1.$

Hence the Laplacian eigenvalues of $M^{t}(n,2m)$ are:
$$5-2cos\frac{2\pi i}{n} \pm\sqrt{3+2cos\frac{2\pi i}{m}
+2cos\frac{2\pi j}{n}+2cos\Big(\frac{2\pi i}{m}+\frac{2\pi
j}{n}\Big)}, \  \ i=0,1,\dots,m-1;j=0,1,\dots,n-1.$$

Similarly, it is not hard to derive the following theorem.
\begin{theorem}\label{4-7}
 For the $M^{t}(n,2m), M^{c}(n,2m)$ and $M^f(n,2m)$ with toroidal,
cylindrical, and free boundary conditions. Then
$$1.~\lim_{m\to\infty} \lim_{n\to
\infty}\frac{LEL\Big(M^t(n,2m)\Big)}{2mn}= \lim_{m\to\infty}
\lim_{n\to \infty}\frac{LEL\Big(M^{c}(n,2m)\Big)}{2mn}=
 \lim_{m\to\infty}
\lim_{n\to \infty}\frac{LEL\Big(M^{f}(n,2m)\Big)}{2mn}\approx
2.1525 .$$\\
$ 2.~LEL \Big(M^t(n,2m)\Big)=LEL \Big(M^c(n,2m)\Big)=LEL
\Big(M^f(n,2m)\Big)\approx 4.3050mn.$
\end{theorem}

{\bf Proof.} By the definition of the Laplacian-energy-like, we
can easily get that
\begin{eqnarray}
\nonumber  LEL\Big(M^t(n,2m)\Big)&=&
\sum_{i=0}^{m-1}\sum_{j=0}^{n-1}
\sqrt{5-2cos\frac{2\pi i}{m} -\sqrt{3+2cos\frac{2\pi i}{m} +2cos\frac{2\pi j}{n}+2cos(\frac{2\pi i}{m}+\frac{2\pi j}{n})}}\\
\nonumber     &&+  \sum_{i=0}^{m-1}\sum_{j=0}^{n-1}
\sqrt{5-2cos\frac{2\pi i}{m} +\sqrt{3+2cos\frac{2\pi i}{m}
+2cos\frac{2\pi j}{n}+2cos(\frac{2\pi i}{m}+\frac{2\pi j}{n})} }.\\
\nonumber
\end{eqnarray}

As an analogue to preceding proof, based on Theorem 3.6, one can
get that
$$\lim_{m\to\infty} \lim_{n\to
\infty}\frac{LEL\Big(M^t(n,2m)\Big)}{2mn}= \lim_{m\to\infty}
\lim_{n\to \infty}\frac{LEL\Big(M^{c}(n,2m)\Big)}{2mn}=
 \lim_{m\to\infty}
\lim_{n\to \infty}\frac{LEL\Big(M^{f}(n,2m)\Big)}{2mn}.$$

It suffices to prove that
\begin{eqnarray}
\nonumber &&\lim_{m\to\infty} \lim_{n\to
\infty}\frac{LEL\Big(M^t(n,2m)\Big)}{2mn}\\ \nonumber &=&
\lim_{m\to\infty} \lim_{n\to
\infty}\frac{1}{2mn}\sum_{i=0}^{m-1}\sum_{j=0}^{n-1} \sqrt{
5-2cos\frac{2\pi i}{n} -\sqrt{3+2cos\frac{2\pi i}{m}
+2cos\frac{2\pi j}{n}+2cos\Big(\frac{2\pi i}{m}+\frac{2\pi
j}{n}\Big)} } \\
 \nonumber & &+\lim_{m\to\infty}
\lim_{n\to
\infty}\frac{1}{2mn}\sum_{i=0}^{m-1}\sum_{j=0}^{n-1}\sqrt{
5-2cos\frac{2\pi i}{n} +\sqrt{3+2cos\frac{2\pi i}{m} +2cos\frac{2\pi
j}{n}+2cos\Big(\frac{2\pi i}{m}+\frac{2\pi
j}{n}\Big)} }\\
   \nonumber &=&\frac{1}{2}\int_{0}^{1}\int_{0}^{1}
\sqrt{ 5-2cos2\pi x -\sqrt{3+2cos 2\pi x
+2cos2\pi y+2cos 2\pi(x+y) }} \ ~dxdy \\
\nonumber & &+\frac{1}{2}\int_{0}^{1}\int_{0}^{1} \sqrt{ 5-2cos2\pi
x +\sqrt{3+2cos 2\pi x
+2cos2\pi y+2cos 2\pi(x+y) }} \ ~dxdy \\
  \nonumber &\approx & 2.1525 .
\end{eqnarray}
The above numerical integration value implies that $M^t(n,2m)$,
$M^{c}(n,2m)$ and $M^{f}(n,2m)$ have the same asymptotic
Laplacian-energy-like, i.e., $LEL \Big(M^t(n,2m)\Big) =LEL
\Big(M^c(n,2m)\Big)=LEL \Big(M^f(n,2m)\Big)\\ \approx 4.3050mn$
 as
$m,n$ tends to infinity. Summing up, we complete the proof. \hfill
$\blacksquare$

\vspace{5pt} \noindent {\bf Remark 4.8}\ In the Figure 2 each kind
of graphs of three kinds of boundary conditions of graphs have the
same the asymptotic Laplacian-energy as $m,n$ tends to infinity.

\section{Concluding remarks}
The calculations of some topological indexes in terms of various
lattices have attracted the attention of many physicists as well
as mathematicians. In this paper, we have deduced the explicit
formulae expressing the Laplacian-energy-like of some lattices
with toroidal, cylindrical, and free boundary conditions, the
explicit asymptotic values of Laplacian-energy-like in these
lattices are obtained via the applications of analysis approach
with the help of calculational software.

 Let $\{G_n\}$ be a
sequence of finite simple graphs with bounded average degree, it
is difficult to calculate its asymptotic Laplacian-energy-like
directly, however, we can find a sequence of graphs $H_n$ with
bounded average degree, which satisfies $\left|V(G_n)\right|$ and
$\left|V(H_n)\right|$ and almost all vertices of $G_n$ and $H_n$
have the same degrees. If we can formulate the asymptotic
Laplacian-energy-like of $H_n$ immediately, then by Theorem 3.6,
$G_n$ and $H_n$ have the same asymptotic Laplacian-energy-like.
Therefore, Theorem 3.6 provides a very effective approach to
handle the asymptotic Laplacian-energy-like of a graph with
bounded average degree. For instance, dealing with the problem of
the asymptotic Laplacian-energy-like of the hexagonal lattice with
the free boundary is not an easy work but we deduced it in a
simple approach.

We can convert some harder problems to easy ones and
simultaneously obtain many results  by utilizing the approach.
Moreover, we showed that the Laplacian-energy-like per vertex of
the many types of lattices is independent of the three boundary
conditions. It is no difficulty to see that the conclusion is true
in general. Actually, the approach can be used widely to formulate
the other topological indexes of various lattices.

 \vspace{5pt} \noindent
{\bf Acknowledgments}\

The work of J. B. Liu is partly supported by the Natural Science
Foundation of Anhui Province of China under Grant No.KJ2013B105;
The work of X. F. Pan is partly supported by the National Science
Foundation of China under Grant Nos.10901001, 11171097, and
11371028.

\end{document}